\documentclass[11pt]{article}
\usepackage{amsthm,amssymb,amsmath,latexsym}
\theoremstyle{definition}
\newtheorem{defin}{Definition}[section]
\theoremstyle{remark}
\newtheorem{example}{Example}

\newtheorem{remar}{Remark}
\theoremstyle{plain}
\newtheorem{thm}[defin]{Theorem}
\newtheorem{prop}[defin]{Proposition}
\newtheorem{lemm}[defin]{Lemma}
\newtheorem{corol}[defin]{Corollary}

\numberwithin{equation}{section}

\begin{document}

\title{Good measures on locally compact Cantor sets}

\author{O.~Karpel\\
Institute for Low Temperature Physics, \\
47 Lenin Avenue, 61103 Kharkov, Ukraine \\
(e-mail: helen.karpel@gmail.com)}

\date{}
\maketitle

\begin{abstract}
We study the set $M(X)$ of full non-atomic Borel (finite or infinite) measures on a non-compact locally compact Cantor set $X$. For an infinite measure $\mu \in M(X)$, the set $\mathfrak{M}_\mu = \{x \in X : \mbox{for any compact open set } U \ni x \mbox{ we have } \mu(U) = \infty \}$ is called defective. We call $\mu$ \textit{non-defective} if $\mu(\mathfrak{M}_\mu) = 0$. The class $M^0(X) \subset M(X)$ consists of probability measures and infinite non-defective measures.
We classify measures $\mu$ from $M^0(X)$ with respect to a homeomorphism. The notions of goodness and compact open values set $S(\mu)$ are defined. A criterion when two good measures from $M^0(X)$ are homeomorphic is given. For any group-like $D \subset [0,1)$ we find a good probability measure $\mu$ on $X$ such that $S(\mu) = D$. For any group-like $D \subset [0,\infty)$ and any locally compact, zero-dimensional, metric space $A$ we find a good non-defective measure $\mu$ on $X$ such that $S(\mu) = D$ and $\mathfrak{M}_\mu$ is homeomorphic to $A$. We consider compactifications $cX$ of $X$ and give a criterion when a good measure $\mu \in M^0(X)$ can be extended to a good measure on $cX$.
\end{abstract}

\section{Introduction}
The problem of classification of Borel finite or infinite measures on topological spaces has a long history. Two measures $\mu$ and $\nu$ defined on Borel subsets of a topological space $X$ are called \textit{homeomorphic} if there exists a self-homeomorphism $h$ of $X$ such that $\mu = \nu\circ h$, i.e. $\mu(E) = \nu(h(E))$ for every Borel subset $E$ of $X$. The topological properties of the space $X$ are important for the classification of measures up to a homeomorphism. For instance, Oxtoby and Ulam \cite{Oxt-Ul} gave a criterion for a Borel probability measure on the finite-dimensional cube to be homeomorphic to the Lebesgue measure. Similar results were obtained for various manifolds (see \cite{Alp-Pr,Oxt-Pr}).

A Cantor set (or Cantor space) is a non-empty zero-dimensional compact perfect metric space. For Cantor sets the situation is much more difficult than for connected spaces. During the last decade, in the papers \cite{Akin3, Austin, S.B.O.K., D-M-Y, Yingst} the Borel probability measures on Cantor sets were studied. In~\cite{K}, infinite Borel measures on Cantor sets were considered. For many applications in dynamical systems the state space is only locally compact. In this paper, we study Borel both finite and infinite measures on non-compact locally compact Cantor sets.

It is possible to construct uncountably many full (the measure of every non-empty open set is positive) non-atomic measures on the Cantor set $X$ which are pairwise non-homeomorphic (see~\cite{Akin1}). This fact is due to the existence of a countable base of clopen subsets of a Cantor set.
The \textit{clopen values set} $S(\mu)$ is the set of finite values of a measure $\mu$ on all clopen subsets of $X$. This set provides an invariant for homeomorphic measures, although it is not a complete invariant.

For the class of the so called \textit{good} probability measures, $S(\mu)$ \textit{is} a complete invariant. By definition, a full non-atomic probability or non-defective measure $\mu$ is good if whenever $U$, $V$ are clopen sets with $\mu(U) < \mu(V)$, there exists a clopen subset $W$ of $V$ such that $\mu(W) = \mu(U)$ (see~\cite{Akin2, K}). Good probability measures are exactly invariant measures of uniquely ergodic minimal homeomorphisms of Cantor sets (see \cite{Akin2}, \cite{GW}).
For an infinite Borel measure $\mu$ on a Cantor set $X$, denote by $\mathfrak{M}_\mu$ the set of all points in $X$ whose clopen neighbourhoods have only infinite measures.
The full non-atomic infinite measures $\mu$ such that $\mu(\mathfrak{M}_\mu) = 0$ are called \textit{non-defective}. These measures arise as ergodic invariant measures for homeomorphisms of a Cantor set and the theory of good probability measures can be extended to the case of non-defective measures (see \cite{K}).

In Section 2, we define a good probability measure and a good non-defective measure on a non-compact locally compact Cantor set $X$ and extend the results concerning good measures on Cantor sets to non-compact locally compact Cantor sets. For a Borel measure $\mu$ on $X$, the set $S(\mu)$ is defined as a set of all finite values of $\mu$ on the compact open sets. The defective set $\mathfrak{M}_\mu$ is the set of all points $x$ in $X$ such that every compact open neighbourhood of $x$ has infinite measure.
We prove the criterion when two good measures on non-compact locally compact Cantor sets are homeomorphic. For every group-like subset $D \subset [0,1)$ we find a good probability measure $\mu$ on a non-compact locally compact Cantor set such that $S(\mu) = D$. For every group-like subset $D \subset [0,\infty)$ and any locally compact, zero-dimensional, metric space $A$ (including $A = \emptyset$) we find a good non-defective measure $\mu$ on a non-compact locally compact Cantor set such that $S(\mu) = D$ and $\mathfrak{M}_\mu$ is homeomorphic to $A$.

In Section 3, compactifications of non-compact locally compact Cantor sets are studied. We investigate whether compactification can be used to classify measures on non-compact locally compact Cantor sets.
We consider only the compactifications which are Cantor sets and extend measure $\mu$ by giving the remainder of compactification a zero measure.
It turns out that in some cases good measure can be extended to a good measure on a Cantor set, while in other cases the extension always produces a measure which is not good. The extensions of a non-good measure are always non-good.
After compactification of a non-compact locally compact Cantor set, new compact open sets are obtained.
We study how the compact open values set changes.
Based on this study, we give a criterion when a good measure on a non-compact locally compact Cantor set stays good after the compactification.

Section 4 illustrates the results of Sections 2 and 3 with the examples. For instance, the Haar measure on the set of $p$-adic numbers and the invariant measure for $(C,F)$-construction are good. We give examples of good ergodic invariant measures on the generating open dense subset of a path space of stationary Bratteli diagrams such that any compactification gives a non-good measure.

\section{Measures on locally compact Cantor sets}

Let $X$ be a non-compact locally compact metrizable space with no isolated points and with a (countable) basis of compact and open sets. Hence $X$ is totally disconnected.
The set $X$ is called a \textit{non-compact locally compact Cantor set}.
Every two non-compact locally compact Cantor sets are homeomorphic (see~\cite{D2}).
Take a countable family of compact open subsets $O_n \subset X$ such that $X = \bigcup_{n=1}^\infty O_n$. Denote $X_1 = O_1$, $X_2 = O_2 \setminus O_1$, $X_3 = O_3 \setminus (O_1 \cup O_2)$,... The subsets $X_n$ are compact, open, pairwise disjoint and $X = \bigcup_{n=1}^\infty X_n$. Since $X$ is non-compact, we may assume without loss in generality that all $X_n$ are nonempty.
Since $X$ has no isolated points, every $X_n$ has the same property. Thus, we represent $X$ as a disjoint union of a countable family of compact Cantor sets $X_n$.

Recall that a Borel measure on a locally compact Cantor space is called \textit{full} if every non-empty open set has a positive measure. It is easy to see that for a non-atomic measure $\mu$ the support of $\mu$ in the induced topology is a locally compact Cantor set. We can consider measures on their supports to obtain full measures. Denote by $M(X)$ the set of full non-atomic Borel measures on $X$.
Then $M(X) = M_f(X) \sqcup M_\infty(X)$, where $M_f(X) = \{\mu \in M(X) : \mu(X) < \infty\}$ and $M_\infty(X) = \{\mu \in M(X) : \mu(X) = \infty\}$. For a measure $\mu \in M_\infty(X)$, denote $\mathfrak{M}_\mu = \{x \in X : \mbox{for any compact and open set } U \ni x \mbox{ we have } \mu(U) = \infty \}$. It will be shown that $\mathfrak{M}_\mu$ is a Borel set. Denote by $M_\infty^{0}(X) = \{\mu \in M_\infty(X) : \mu(\mathfrak{M}_\mu) = 0\}$. Let $M^0(X) = M_f(X) \sqcup M_\infty^{0}(X)$. Throughout the paper we will consider only measures from $M^0(X)$. We normalize the measures from $M_f(X)$ so that $\mu(X) = 1$ for any $\mu \in M_f(X)$.

Recall that $\mu \in M^0(X)$ is \textit{locally finite} if every point of X has a neighbourhood of finite measure. The properties of measures from the class $M^0(X)$ are collected in the following proposition.

\begin{prop}\label{basic_prop}
Let $\mu \in M^0(X)$. Then

\noindent(1) The measure $\mu$ is locally finite if and only if $\mathfrak{M}_\mu = \emptyset$,

\noindent(2) The set $X \setminus \mathfrak{M}_\mu$ is open. The set $\mathfrak{M}_\mu$ is $F_\sigma$.

\noindent(3) For any compact open set $U$ with $\mu(U) = \infty$ and any $a > 0$ there exists a compact open subset $V \subset U$ such that $a \leq \mu(V) < \infty$.

\noindent(4)  The set $\mathfrak{M}_\mu$ is nowhere dense.

\noindent(5) $X = \bigsqcup_{i = 1}^\infty V_i \bigsqcup \mathfrak{M}_\mu$, where each $V_i$ is a compact open set of finite measure and $\mathfrak{M}_\mu$ is a nowhere dense $F_\sigma$ and has zero measure. The measure $\mu$ is $\sigma$-finite.

\noindent(6) $\mu$ is uniquely determined by its values on the algebra of compact open sets.

\end{prop}

\noindent \textbf{Proof.} (1) The condition $\mathfrak{M}_\mu = \emptyset$ means that every point $x \in X$ has a compact open neighbourhood of finite measure.
Hence $\mu$ is locally finite and vise versa.

(2) We have $X \setminus \mathfrak{M}_\mu = \{x \in X : \mbox{ there exists a compact open set } U_x \ni x \mbox{ such that } \mu(U_x) < \infty \}$. Then for every point $x \in X \setminus \mathfrak{M}_\mu$ we have $U_x \subset X \setminus \mathfrak{M}_\mu$. Hence $X \setminus \mathfrak{M}_\mu$ is open. Therefore, for every $n \in \mathbb{N}$ the set $X_n \setminus \mathfrak{M}_\mu$ is open and $X_n \cap \mathfrak{M}_\mu$ is closed. Then $\mathfrak{M}_\mu = \bigsqcup_{n \in \mathbb{N}} \; (X_n \cap \mathfrak{M}_\mu)$ is $F_\sigma$ set.

(3) Let $U$ be a non-empty compact open subset of $X$ such that $\mu(U) = \infty$. Since $\mu \in M^{0}(X)$, we have $\mu(U) = \mu(U \setminus \mathfrak{M}_{\mu})$. Since $U$ is open, the set $U \setminus \mathfrak{M}_{\mu} = U \cap (X \setminus \mathfrak{M}_{\mu})$ is open.
There are only countably many compact open subsets in $X$, hence the open set $U \setminus \mathfrak{M}_{\mu}$ can be represented as a disjoint union of compact open subsets $\{U_i\}_{i \in \mathbb{N}}$ of finite measure.
We have $\mu(U) = \sum_{i=0}^{\infty} \mu(U_i) = \infty$, hence for every $a \in \mathbb{R}$ there is a compact open subset $V = \bigsqcup_{i=0}^{N}U_i$ such that $a \leq \mu(V) < \infty$.

(4) Let $U$ be a compact open subset of $X$. It suffices to show that there exists a non-empty compact open subset $V \subset U$ such that $V \cap \mathfrak{M}_{\mu} = \emptyset$. If $\mu(U) < \infty$ then $U \cap \mathfrak{M}_{\mu} = \emptyset$. Otherwise, by (3), there exists a compact open subset $V \subset U$ such that $0 < \mu(V) < \infty$. Obviously, $V \cap \mathfrak{M}_{\mu} = \emptyset$.

(5) follows from the proof of (3).

(6) follows from (5). 
$\blacksquare$

\medskip
For a measure $\mu \in M^0(X)$ define the \textit{compact open values set} as the set of all finite values of the measure $\mu$ on the compact open sets:
$$
S(\mu) = \{\mu(U):\,U\mbox{ is compact open in } X \mbox{ and } \mu(U) < \infty\}.
$$
For each measure $\mu \in M^0(X)$, the set $S(\mu)$ is a countable dense subset of the interval $[0, \mu(X))$. Indeed, the set $S(\mu)$ is dense in $[0, \mu(V)]$ for every compact open set $V$ of finite measure (see~\cite{Akin1}). By Proposition~\ref{basic_prop}, $S(\mu)$ is dense in $[0, \mu(X))$.

Let $X_{1}$, $X_{2}$ be two non-compact locally compact Cantor sets.
It is said that measures $\mu_{1} \in M(X_{1})$ and $\mu_{2} \in M(X_{2})$ are \textit{homeomorphic} if there exists a homeomorphism $h \colon X_{1} \rightarrow X_{2}$ such that $\mu_{1}(E) = \mu_{2}(h(E))$ for every Borel subset $E \subset X_1$.
Clearly, $S(\mu_{1}) = S(\mu_{2})$ for any homeomorphic measures $\mu_1$ and $\mu_2$.
We call two Borel infinite measures $\mu_1 \in M^0_\infty(X_{1})$ and $\mu_2 \in M^0_\infty(X_{2})$ \textit{weakly homeomorphic} if there exists a homeomorphism $h \colon X_{1} \rightarrow X_{2}$ and a constant $C>0$ such that $\mu_{1}(E) = C \mu_{2}(h(E))$ for every Borel subset $E \subset X_1$. Then $S(\mu_{1}) = C S(\mu_{2})$.

Let $D$ be a dense countable subset of the interval $[0,a)$ where $a \in (0, \infty]$.
Then $D$ is called \textit{group-like} if there exists an additive subgroup $G$ of $\mathbb{R}$ such that $D = G \cap [0, a)$. It is easy to see that $D$ is group-like if and only if for any $\alpha, \beta \in D$ such that $\alpha \leq \beta$ we have $\beta - \alpha \in D$ (see \cite{Akin2, K}). 

\begin{defin}
Let $X$ be a locally compact Cantor space (either compact or non-compact) and $\mu\in M^0(X)$.
A compact open subset $V$ of $X$ is called \textit{good} for $\mu$ (or just good when the measure is understood) if for every compact open subset $U$ of $X$ with $\mu(U) < \mu(V)$, there exists a compact open set $W$ such that $W \subset V$ and $\mu(W) = \mu(U)$. A measure $\mu$ is called \textit{good} if every compact open subset of $X$ is good for $\mu$.
\end{defin}

If $\mu \in M^0(X)$ is a good measure and $\nu \in M^0(X)$ is (weakly) homeomorphic to $\mu$ then, obviously, $\nu$ is good.
It is easy to see that in the case of compact Cantor set the definition of a good measure coincides with the one given in~\cite{Akin2}. For a compact open subset $U \subset X$ let $\mu|_U$ be the restriction of the measure $\mu$ to the Cantor space $U$. Then the set $U$ is good if and only if $S(\mu|_U) = S(\mu|_X) \cap [0, \mu(U)]$. Denote by $H_{\mu}(X)$ the group of all homeomorphisms of a space $X$ preserving the measure $\mu$. The action of $H_{\mu}(X)$ on $X$ is called \textit{transitive} if for every $x_{1}, x_{2} \in X$ there exists $h \in H_{\mu}(X)$ such that $h(x_{1}) =  x_{2}$. The action is called \textit{topologically transitive} if there exists a dense orbit, i.e. there is $x \in X$ such that the set $O(x) =  \{h(x) : h \in H_\mu(X)\}$ is dense in $X$.

We extend naturally the notion of partition basis introduced in~\cite{Akin3}.
A \textit{partition basis} $\mathcal{B}$ for a non-compact locally compact Cantor set $X$ is a collection of compact open subsets of $X$ such that every non-empty compact open subset of $X$ can be partitioned by elements of $\mathcal{B}$.

The properties of good measures on non-compact locally compact Cantor sets are gathered in the following proposition. The proofs for the measures on compact Cantor spaces can be found in~\cite{Akin2,Akin3, K}.

\begin{prop}\label{many}
Let $X$ be a locally compact Cantor space (either compact or non-compact). Let $\mu \in M^0(X)$. Then

(a) If $\mu$ is good and $C > 0$ then $C \mu$ is good and $S(C \mu) = C S(\mu)$.

(b) If $\mu$ is good and $U$ is a non-empty compact open subset of $X$ then the measure $\mu|_U$ is good and $S(\mu|_U) = S(\mu) \cap [0, \mu(U)]$.

(c) $\mu$ is good if and only if every compact open subset of finite measure is good.

(d) $\mu$ is good if and only if for every non-empty compact open subset $U$ of finite measure, the measure $\mu|_U$ is good.

(e) If $\mu$ is good then $S(\mu)$ is group-like.

(f) If a compact open set $U$ admits a partition by good compact open subsets then $U$ is good.

(g) The measure $\mu$ is good if and only if there exists a partition basis $\mathcal B$ consisting of compact open sets which are good for $\mu$.

(h) If $\mu$ is good, then the group $H_\mu(X)$ acts transitively on $X \setminus \mathfrak{M_\mu}$. In particular, the group $H_\mu(X)$ acts topologically transitively on $X$.

(i) If $\mu$ is a good measure on $X$ and $\nu$ is the counting measure on $\{1,2,...,n\}$ then $\mu \times \nu$ is a good measure on $X \times \{1,2,...,n\}$.

\end{prop}

\noindent \textbf{Proof}. (a), (b) are clear.

(c) Suppose that every compact open subset of finite measure is good. Let $V$ be any compact open set with $\mu(V) = \infty$ and $U$ be a compact open set with $\mu(U) < \infty$. By Proposition~\ref{basic_prop}, there exists a compact open subset $W \subset V$ such that $\mu(U) \leq \mu(W) < \infty$. By assumption, $W$ is good. Hence there exists a compact open set $W_1\subset W$ with $\mu(W_1) = \mu(U)$ and $V$ is good.

(d) Suppose that for every non-empty compact open subset $U$ of finite measure, the measure $\mu|_U$ is good. We prove that every compact open subset of finite measure is good, then use (c). Let $U$, $V$ be compact open sets with $0 < \mu(U) < \mu(V) < \infty$. Set $W = U \cup V$. Then $W$ is a compact open set of finite measure. Since $\mu|_W$ is good, there exists $W_1 \subset V$ such that $\mu(W_1) = \mu(U)$.

(e) If $\mu$ is good then for any $\alpha, \beta \in S(\mu)$ such that $\beta - \alpha \geq 0$, we have $\beta - \alpha \in S(\mu)$. Hence $S(\mu)$ is group-like. 

(f) See~\cite{Akin3} for the case of finite measure and~\cite{K} for infinite measure.

(g) If there exists a partition basis $\mathcal B$ consisting of compact open sets which are good for $\mu$, then, by (f), every compact open set is good.

(h) For any $x, y \in X \setminus \mathfrak{M_\mu}$ there exists a compact open set $U$ of finite measure such that $x, y \in U$. By (d), the measure $\mu|_U$ is a good finite measure on a Cantor space $U$. By Theorem 2.13 in~\cite{Akin2}, there exists a homeomorphism $h \colon U \rightarrow U$ which preserves $\mu$ and $h(x) = y$. Define $h_1 \in H_\mu(X)$ to be $h$ on $U$ and the identity on $X \setminus U$. For every $x \in X\setminus\mathfrak{M}_\mu$ we have $O(x) = X \setminus\mathfrak{M}_\mu$. By Proposition~\ref{basic_prop}, the set $X \setminus\mathfrak{M}_\mu$ is dense in $X$. Hence $H_\mu(X)$ acts topologically transitively on $X$.

(i) The rectangular compact open sets $U \times \{z\}$, where $U$ is compact open in $X$ and $z \in \{1,2,...,n\}$, form a partition basis for $X \times \{1,2,...,n\}$. Since $\mu \times \nu (U \times \{z\}) = \mu(U)$, these sets are good. The measure $\mu$ is good by (g).
$\blacksquare$

\medskip
For $G$ an additive subgroup of $\mathbb{R}$ we call a positive real number $\delta$ a \textit{divisor} of $G$ if $\delta G = G$. The set of all divisors of $G$ is called $Div(G)$. By a full measure on a discrete countable topological space $Y$ we mean a measure $\nu$ such that $0 < \nu(\{y\}) < \infty$ for every $y \in Y$. We will use the following theorem for $Y = \mathbb{Z}$, but the proof stays correct for any discrete countable topological space $Y$.

\begin{thm}\label{good_product}
Let $\mu$ be a good measure on a non-compact locally compact Cantor space $X$. Let $\nu$ be a full measure on $\mathbb{Z}$, where $\mathbb{Z}$ is endowed with discrete topology. Let $G$ be an additive subgroup of $\mathbb{R}$ generated by $S(\mu)$.
Then $\mu \times \nu$ is good on $X \times \mathbb{Z}$ if and only if there exists $C > 0$ such that $\nu(\{i\}) \in C \cdot Div(G)$ for every $i \in \mathbb{Z}$.
\end{thm}

\noindent \textbf{Proof.} Lets prove the ``if'' part. Suppose $\mu$ is good on $X$ and $\nu(\{i\}) \in C \cdot Div(G)$ for some $C > 0$ and every $i \in \mathbb{Z}$.
By Proposition~\ref{many} (g), it suffices to prove that a compact open set of the form $U \times \{i\}$ is good for any compact open $U \subset X$ and any $i \in \mathbb{Z}$. Thus, it suffices to show that $S(\mu \times \nu|_{U \times \{i\}}) = S(\mu \times \nu|_{X \times \mathbb{Z}}) \cap [0, \mu \times \nu (U \times \{i\})]$. The inclusion $S(\mu \times \nu|_{U \times \{i\}}) \subset S(\mu \times \nu|_{X \times \mathbb{Z}}) \cap [0, \mu \times \nu (U \times \{i\})]$ is always true, hence we need to prove the inverse inclusion.
We have $S(\mu \times \nu|_{U \times\{i\}}) = \nu(\{i\}) S(\mu|_U) = C \delta S(\mu|_U)$ for some $\delta \in Div(G)$. Since $\mu$ is good on $X$, we obtain $S(\mu|_U) = G \cap [0, \mu(U)]$. Hence $S(\mu \times \nu|_{U \times\{i\}}) = C G \cap [0, C \delta \mu(U)] = C G \cap [0, \mu \times \nu (U \times \{i\})]$. Note that $C \delta \mu(U) \in CG$ because $\delta \in Div(G)$.
Therefore, it suffices to prove that $S(\mu \times \nu|_{X \times \mathbb{Z}}) \subset C G$. The set $S(\mu \times \nu|_{X \times \mathbb{Z}})$ consists of all finite sums $\sum_{i,j} \mu(U_i) \nu(\{j\})$, where each $U_i$ is a compact open set in $X$ and $j \in \mathbb{Z}$. We have $\sum_{i,j} \mu(U_i) \nu(\{j\}) = \sum_{i,j} \mu(U_i) C \delta_j \subset CG$, here $\delta_i \in Div(G)$. Hence $S(\mu \times \nu|_{U \times\{i\}}) \supset S(\mu \times \nu|_{X \times \mathbb{Z}}) \cap [0, \mu \times \nu (U \times \{i\})]$ and $U \times \{i\}$ is good.

Now we prove the ``only if part''. Suppose that $\mu \times \nu$ is good on $X \times \mathbb{Z}$. Then for any $i \in \mathbb{Z}$ we have $S(\mu \times \nu|_{X \times\{i\}}) = S(\mu \times \nu|_{X \times \mathbb{Z}}) \cap [0, \mu \times \nu(X \times\{i\})]$. Note that $S(\mu \times \nu|_{X \times\{i\}} = \nu(\{i\}) S(\mu|_X)$. Denote by $\widetilde{G}$ the additive subgroup of $\mathbb{R}$ generated by $S(\mu \times \nu|_{X \times \mathbb{Z}})$. Let $\alpha = \nu(\{i\})$. Then $\alpha G = \widetilde{G}$. Let $j \in \mathbb{Z}$ and $\beta = \nu(\{j\})$. By the same arguments, we have $\beta G = \widetilde{G}$. Then $\frac{\alpha}{\beta} \in Div(G)$. Indeed, $\frac{\alpha}{\beta} G = \frac{1}{\beta} \widetilde{G} = G$. Hence $\alpha = \beta \delta$, where $\delta \in Div(G)$. Set $C = \nu(\{j\})$. Then for every $i \in \mathbb{Z}$ we have $\nu(\{i\}) = C \delta_i$ where $\delta_i = \frac{\nu(\{i\})}{\nu(\{j\})} \in Div (G)$.
$\blacksquare$

\begin{thm}
Let $X$, $Y$ be non-compact locally compact Cantor sets.
If $\mu \in M^0(X)$, $\nu \in M^0(Y)$ are good measures, then the product $\mu \times \nu$ is a good measure on $X \times Y$ and
$$
S(\mu \times \nu) = \left\{\sum_{i=0}^N \alpha_i \cdot \beta_i : \alpha_i \in S(\mu), \beta_i \in S(\nu), N \in \mathbb{N}\right\} \cap [0, \mu(X)\times \nu(Y)).
$$
\end{thm}

\noindent \textbf{Proof.} Let $X = \bigsqcup_{m = 1}^{\infty} X_n$ and $Y = \bigsqcup_{n = 1}^{\infty} Y_n$, where each $X_n$, $Y_n$ is a Cantor set. Then $X \times Y = \bigsqcup_{m, n = 1}^{\infty} X_m \times Y_n$ and $\mu \times \nu |_{X_m \times Y_n} = \mu|_{X_n} \times \nu|_{Y_n}$. Since $\mu|_{X_n}$ and $\nu|_{Y_n}$ are good finite or non-defective measures on a Cantor set, the measure $\mu \times \nu |_{X_m \times Y_n}$ is good by Theorem 2.8 (\cite{Akin3}), Theorem 2.10 (\cite{K}). By Proposition~\ref{many}, $\mu \times \nu$ is good on $X \times Y$. $\blacksquare$

\begin{thm}\label{krit_homeo_good}
Let $X$, $Y$ be non-compact locally compact Cantor spaces. Let $\mu \in M^0(X)$ and $\nu \in M^0(Y)$ be good measures. Let $S(\mu) = S(\nu)$. Let $\mathfrak{M}$ be the defective set for $\mu$ and $\mathfrak{N}$ be the defective set for $\nu$. Assume that there is a homeomorphism $h \colon \mathfrak{M} \rightarrow \mathfrak{N}$ where the sets $\mathfrak{M}$ and $\mathfrak{N}$ are endowed with the induced topologies. Then there exists a homeomorphism $\widetilde{h} \colon X \rightarrow Y$ which extends $h$ such that $\mu = \nu \circ \widetilde{h}$.

Conversely, if $\mu \in M^0(X)$ and $\nu \in M^0(Y)$ are good homeomorphic measures then $S(\mu) = S(\nu)$ and there is a homeomorphism $h \colon \mathfrak{M} \rightarrow \mathfrak{N}$.
\end{thm}

\noindent \textbf{Proof.} The second part of the Theorem is clear. We prove the first part.
Let $X = \bigsqcup_{i=1}^\infty X_i$ and $Y = \bigsqcup_{j=1}^\infty Y_j$ where $X_i$, $Y_j$ are compact Cantor spaces.

First, consider the case when $\mathfrak{M} = \mathfrak{N} = \emptyset$, i.e. the measures $\mu$, $\nu$ are either finite of infinite locally finite measures. Since $S(\mu) = S(\nu)$, we have $\mu(X_1) \in S(\nu)$. There exists $n \in \mathbb{N}$ such that $\nu(\bigsqcup_{j=1}^{n-1} Y_j) \leq \mu(X_1) < \nu(\bigsqcup_{j=1}^{n} Y_j)$. Since $S(\nu)$ is group-like, we see that $\mu(X_1) - \nu(\bigsqcup_{j=1}^{n-1} Y_j) \in S(\nu)$. Since $\nu$ is good, there exists a compact open subset $W \subset Y_n$ such that $\nu(W) = \mu(X_1) - \nu(\bigsqcup_{j=1}^{n-1} Y_j)$. Hence $Z = \bigsqcup_{j=1}^{n-1} Y_j \sqcup W$ is a compact Cantor set and $\mu(X_1) = \nu(Z)$. By Theorem 2.9 (\cite{Akin2}), there exists a homeomorphism $h_1 \colon X_1\rightarrow Z$ such that $\mu |_{X_1} = \nu |_{Z} \circ h_1$. Set $\widetilde{h}|_{X_1} = h_1$. Consider $(Y_n \setminus W) \bigsqcup_{j=n+1}^\infty Y_j$ instead of $Y$ and $\bigsqcup_{i=2}^\infty X_i$ instead of $X$. Reverse the roles of $X$ and $Y$. Proceed in the same way using $Y_n \setminus W$ instead of $X_1$. Thus, we obtain countably many homeomorphisms $\{h_i\}_{i=1}^\infty$. Given $x \in X$, set $\widetilde{h}(x) = h_i(x)$ for the corresponding $h_i$. Then $\widetilde{h} \colon X \rightarrow Y$ is a homeomorphism which maps $\mu$ into $\nu$.

Now, let $\mathfrak{M} \neq \emptyset$. If $\mu(X_1) < \infty$, we proceed as in the previous case. If $\mu(X_1) = \infty$ then $X_1 \cap \mathfrak{M} \neq \emptyset$. Then $h(X_1 \cap \mathfrak{M})$ is a compact open subset of $\mathfrak{N}$ in the induced topology. Hence there exists a compact open set $W \subset Y$ such that $W \cap \mathfrak{N} = h(X_1 \cap \mathfrak{M})$. Then, by Theorem 2.11 (\cite{K}), the sets $X_1$ and $W$ are homeomorphic via measure preserving homeomorphism $h_1$ and $h_1|_{X_1 \cap \mathfrak{M}} = h$. Since $W$ is compact, there exists $N$ such that $W \subset \bigsqcup_{n=1}^N Y_n$. Reverse the roles of $X$ and $Y$ and consider $\bigsqcup_{n=1}^N Y_n \setminus W$ instead of $X_1$. $\blacksquare$

The corollary for weakly homeomorphic measures follows:
\begin{corol}\label{krit_weak_homeo_good}
Let $\mu \in M_{\infty}^0(X)$ and $\nu \in M_{\infty}^0(Y)$ be good infinite measures on non-compact locally compact Cantor sets $X$ and $Y$. Let $\mathfrak{M}$ be the defective set for $\mu$ and $\mathfrak{N}$ be the defective set for $\nu$. Then $\mu$ is weakly homeomorphic to $\nu$ if and only if the following conditions hold:

(1) There exists $c > 0$ such that $S(\mu) = c S(\nu)$,

(2) There exists a homeomorphism $h \colon \mathfrak{M} \rightarrow \mathfrak{N}$ where the sets $\mathfrak{M}$ and $\mathfrak{N}$ are endowed with the induced topologies.

\end{corol}

\begin{remar}
Let $\mu \in M^0_{\infty}(X)$ be a good measure on a non-compact locally compact Cantor set $X$ and $V$ be any compact open subset of $X$ with $\mu(V) < \infty$. Then $\mu$ on $X$ is homeomorphic to $\mu$ on $X \setminus V$.
Let $S(\mu) = G \cap [0, \infty)$. Then $\mu$ is homeomorphic to $c \mu$ if and only if $c \in Div(G)$.
\end{remar}

\begin{corol}
Let $\mu$ be a good finite or non-defective measure on a non-compact locally compact Cantor set $X$. Let $U$, $V$ be two compact open subsets of $X$ such that $\mu(U) = \nu(V) < \infty$. Then there is $h \in H_\mu(X)$ such that $h(U) = V$.
\end{corol}

\noindent \textbf{Proof.}
Set $Y = U \cup V$. Then $Y$ is a Cantor set with $\mu(Y) < \infty$. By Proposition 2.11 in~\cite{Akin2}, there exists a self-homeomorphism $h$ of $Y$ such that $h(U) = V$ and $h$ preserves $\mu$. Set $h$ to be identity on $X \setminus Y$. $\blacksquare$

\begin{corol}
Let $\mu$ and $\nu$ be good non-defective measures on non-compact locally compact Cantor sets $X$ and $Y$. Let $\mathfrak{M}$ be the defective set for $\mu$ and $\mathfrak{N}$ be the defective set for $\nu$. If there exist compact open sets $U \subset X$ and $V \subset Y$ such that $\mu(U) = \nu(V) < \infty$ and $\mu|U$ is homeomorphic to $\nu|V$, then $\mu$ is homeomorphic to $\nu$ if and only if $\mathfrak{M}$ and $\mathfrak{N}$ (with the induced topologies) are homeomorphic.
\end{corol}

\noindent \textbf{Proof.} Let $\gamma = \mu(U) = \nu (V)$. Since $\mu|U$ is homeomorphic to $\nu|V$, we have $S(\mu|U) = S(\nu|V)$. Since $\mu$ and $\nu$ are good, we have $S(\mu) \cap [0, \gamma] = S(\nu) \cap [0, \gamma]$ by Proposition~\ref{many}. Since $S(\mu)$ and $S(\nu)$ are group-like, we obtain $S(\mu) = S(\nu)$. $\blacksquare$

\begin{thm}\label{goodSmu}
Let $\mu \in M^{0}(X)$ be a good measure on a non-compact locally compact Cantor set $X$. Then the compact open values set $S(\mu)$ is group-like and the defective set $\mathfrak{M}_\mu$ is a locally compact, zero-dimensional, metric space (including $\emptyset$).

Conversely, for every countable dense group-like subset $D$ of $[0, 1)$, there is a good probability measure $\mu$ on a non-compact locally compact Cantor set such that $S(\mu) = D$. For every countable dense group-like subset $D$ of $[0, \infty)$ and any locally compact, zero-dimensional, metric space $A$ (including $A = \emptyset$) there is a good non-defective measure $\mu$ on a non-compact locally compact Cantor set such that $S(\mu) = D$ and $\mathfrak{M}_\mu$ is homeomorphic to $A$.
\end{thm}

\noindent \textbf{Proof}. The first part of the theorem follows from Propositions~\ref{basic_prop},~\ref{many}.

We prove the second part. First, consider the case of finite measure. Let $D \subset [0,1)$ be a countable dense group-like subset. Then there exist a strictly increasing sequence $\{\gamma_n\}_{n=1}^{\infty} \subset D$ such that $\lim_{n \rightarrow \infty} \gamma_n = 1$. For $n = 1,2,...$ set $\delta_n = \gamma_{n} - \gamma_{n-1}$.
Denote by $S_n = D \cap [0, \delta_n]$. Then $D_n = \frac{1}{\delta_n} (D \cap [0,\delta_n])$ is a group-like subset of $[0,1]$ with $1 \in D_n$. In~\cite{Akin2}, it was proved that there exists a good probability measure $\mu_n$ on a Cantor set $X_n$ such that $S(\mu_n|_{X_n}) = D_n$. The measure $\nu_n = \delta_n \mu_n$ is a good finite measure on $X_n$ with $S(\nu_n|_{X_n}) = D \cap [0,\delta_n]$. Set $X = \bigsqcup_{n=1}^{\infty} X_n$ and let $\mu|_{X_n} = \nu_n$. Then $\mu$ is a good probability measure on a non-compact locally compact Cantor space $X$ and $S(\mu|_{X}) = D$.

Now consider the case of infinite measure. Let $\gamma \in D$. Since $D \subset [0, \infty)$ is group-like, we see that $\frac{1}{\gamma}D \cap [0,1]$ is a group-like subset of $[0,1]$. In~\cite{Akin2} it was proved that there exists a good probability measure $\mu_1$ on a Cantor space $Y$ with $S(\mu_1) = \frac{1}{\gamma}D \cap [0,1]$. Set $\mu = \gamma \mu_1$. Then $\mu$ is a good finite measure on $Y$ and $S(\mu) = D \cap [0,\gamma]$. Endow the set $\mathbb{Z}$ with discrete topology. Let $\nu$ be a counting measure on $\mathbb{Z}$. Set $X = Y \times \mathbb{Z}$ and $\widetilde{\mu} = \mu \times \nu$. Then, by Theorem~\ref{good_product}, $\widetilde{\mu}$ is good with $S(\widetilde{\mu}) = D$ and $\mathfrak{M}_{\widetilde{\mu}} = \emptyset$.

Suppose $A$ is a non-empty compact zero-dimensional, metric space. Then, by Theorem 2.15 (\cite{K}), there exists a good non-defective measure $\mu$ on a Cantor space $Y$ such that $S(\mu) = D$ and $\mathfrak{M}_\mu$ is homeomorphic to $A$. By the above, there exists a good locally finite measure $\nu$ on a non-compact locally compact set $X$ with $S(\nu) = D$ and $\mathfrak{M}_\nu = \emptyset$. Set $Z = Y \sqcup X$ and $\widetilde{\mu}|_{Y} = \mu$, $\widetilde{\mu}|_{X} = \nu$. Then $\widetilde{\mu}$ is good on a non-compact locally compact Cantor set $Z$ with $S(\widetilde{\mu}) = D$ and $\mathfrak{M}_{\widetilde{\mu}}$ is homeomorphic to $A$.

Suppose that $A$ is a non-empty, non-compact, locally compact, zero-dimensional metric space. Then $A = \bigsqcup_{n=1}^{\infty} A_n$ where each $A_n$ is a non-empty, compact, zero-dimensional metric space. By Theorem 2.15 (\cite{K}), for every $n = 1,2,...$ there exists a good non-defective measure $\mu_n$ on a Cantor set $Y_n$ such that $S(\mu_n) = D$ and $\mathfrak{M}_{\mu_n}$ is homeomorphic to $A_n$. Set $X = \bigsqcup_{n=1}^{\infty} Y_n$ and $\mu|_{Y_n} = \mu_n$. Then $\mu$ is good on a non-compact locally compact Cantor set $X$ with $S(\mu) = D$ and $\mathfrak{M}_{\widetilde{\mu}}$ is homeomorphic to $A$. $\blacksquare$

\begin{corol}\label{invarhomeo}
Let $D$ be a countable dense group-like subset of $[0, \infty)$. Then there exists an aperiodic homeomorphism of a non-compact locally compact Cantor set with good non-defective invariant measure $\widetilde{\mu}$ such that $S(\widetilde{\mu}) = D$.
\end{corol}

\noindent \textbf{Proof}. We use the construction similar to the one in the proof of Theorem~\ref{goodSmu}. Let $\mu$ be a good measure on a Cantor set $Y$ with $S(\mu) = D \cap [0,\gamma]$ for some $\gamma \in D$. Let $\nu$ be a counting measure on $\mathbb{Z}$.
Set $\widetilde{\mu} = \mu \times \nu$ on $X = Y \times \mathbb{Z}$. Then $\widetilde{\mu}$ is a good non-defective measure on a non-compact locally compact Cantor set $X$ with $S(\widetilde{\mu}) = D$.
Since the measure $\mu$ is a good finite measure on $Y$, there exists a minimal homeomorphism $T \colon Y \rightarrow Y$ such that $\mu$ is invariant for $T$ (see~\cite{Akin2}). Let $T_1(x,n) = (Tx, n+1)$. Then $T_1$ is aperiodic homeomorphism of $X$. The measure $\widetilde{\mu}$ is invariant for $T_1$. $\blacksquare$

\begin{remar}
The measure $\widetilde{\mu}$ built in Corollary~\ref{invarhomeo} is invariant for any skew-product with the base $(Y,T)$ and cocycle acting on $\mathbb{Z}$.
\end{remar}

\begin{thm}
Let $X$ be a non-compact locally compact Cantor set.
Then there exist continuum distinct classes of homeomorphic good measures in $M_f(X)$.
There also exist continuum distinct classes of weakly homeomorphic good measures in $M_\infty^0(X)$.
\end{thm}

\noindent \textbf{Proof.} There exist uncountably many distinct group-like subsets $\{D_\alpha\}_{\alpha \in \Lambda}$ of $[0,1]$. By Theorem~\ref{goodSmu}, for each $D_\alpha$ there exists a good probability measure $\mu_\alpha$ on $X$ such that $S(\mu_\alpha) = D_\alpha$. By Theorem~\ref{krit_homeo_good}, the measures $\{\mu_\alpha\}_{\alpha \in \Lambda}$ are pairwise non-homeomorphic.

Let $Y$ be a compact Cantor set. Let $\mu$ be a non-defective measure on $Y$. Denote by $[\mu]$ the class of weak equivalence of $\mu$ in the set of all non-defective measures on $Y$. There exist continuum distinct classes $[\mu_\alpha]$ of weakly homeomorphic good non-defective measures on a Cantor set $Y$ (see Theorem 2.18 in~\cite{K}). Moreover, if there exists $C>0$ such that $G(S(\mu_\alpha)) = C G(S(\mu_\beta))$ then $\mu_\beta \in [\mu_\alpha]$.
Let $\nu$ be a counting measure on $\mathbb{Z}$. Then, by Theorem~\ref{good_product}, $\mu_\alpha \times \nu$ is a good measure on a non-compact locally compact Cantor set $Y \times \mathbb{Z}$ and $G(S(\mu_\alpha \times \nu)) = G(S(\mu_\alpha))$. Hence, by Corollary~\ref{krit_weak_homeo_good}, the measures $\mu_\alpha \times \nu$ and $\mu_\beta \times \nu$ are weakly homeomorphic if and only if $\mu_\beta \in [\mu_\alpha]$.
$\blacksquare$

\begin{prop}
If $\mu$ is Haar measure for some topological group structure on a non-compact locally compact Cantor space $X$ then $\mu$ is a good measure on $X$.
\end{prop}

\noindent \textbf{Proof.} The ball $B$ centered at the identity in the invariant ultrametric is a compact open subgroup of $X$. Since $\mu$ is translation-invariant, by Proposition~\ref{many}, it suffices to show that $\mu|_B$ is good for every such ball $B$. Since the restriction of $\mu$ on $B$ is a Haar measure on a compact Cantor space, $\mu|_B$ is good by Proposition 2.4 in~\cite{Akin3}. $\blacksquare$

\section{From measures on non-compact spaces to measures on compact spaces and back again}

Let $X$ be a non-compact locally compact Cantor space. A \textit{compactification} of $X$ is a pair $(Y,c)$ where $Y$ is a compact space and $c \colon X \rightarrow Y$ is a homeomorphic embedding of $X$ into $Y$ (i.e. $c \colon X \rightarrow c(X)$ is a homeomorphism) such that $\overline{c(X)} = Y$, where $\overline{c(X)}$ is the closure of $c(X)$. In the paper, by compactification we will mean not only a pair $(Y,c)$ but also the compact space $Y$ in which $X$ can be embedded as a dense subset. We will denote the compactifications of a space $X$ by symbols $cX$, $\omega X$, etc., where $c$, $\omega$ are the corresponding homeomorphic embeddings.

Let $\mu \in M^0(X)$.
We will consider only such compactifications $cX$ that $cX$ is a Cantor set. Since $c$ is a homeomorphism, the measure $\mu$ on $X$ passes to a homeomorphic measure on $c(X)$. Since we are interested in the classification of measures up to homeomorphisms, we can identify the set $c(X)$ with $X$. Hence $X$ can be considered as an open dense subset of $cX$. The set $cX \setminus X$ is called the \textit{remainder} of compactification. As far as $X$ is locally compact, the remainder $cX \setminus X$ is closed in $cX$ for every compactification $cX$ (see~\cite{E}). Since $\overline{X} = cX$, the set $cX \setminus X$ is a closed nowhere dense subset of $cX$.

Compactifications $c_1X$ and $c_2X$ of a space $X$ are \textit{equivalent} if there exists a homeomorphism $f \colon c_1X \rightarrow c_2X$ such that $fc_1(x) = c_2(x)$ for every $x \in X$. We shall identify equivalent compactifications.
For any space $X$ one can consider the family $\mathcal{C}(X)$ of all compactifications of $X$. The order relation on $\mathcal{C}(X)$ is defined as follows: $c_2X \leq c_1X$ if there exists a continuous map $f \colon c_1X \rightarrow c_2X$ such that $fc_1 = c_2$. Then we have $f(c_1(X)) = c_2(X)$ and $f(c_1X \setminus c_1(X)) = c_2X \setminus c_2(X)$.

\begin{thm}[The Alexandroff compactification theorem]\label{Alex}
Every non-compact locally compact space $X$ has a compactification $\omega X$ with one-point remainder. This compactification is the smallest element in the set of all compactifications $\mathcal{C}(X)$ with respect to the order $\leq$.
\end{thm}

The topology on $\omega X$ is defined as follows. Denote by $\{\infty\}$ the point $\omega X \setminus X$. Open sets in $\omega X$ are the sets of the form $\{\infty\} \cup (X \setminus F)$, where $F$ is a compact subspace of $X$, together with all sets that are open in $X$.

For any Borel measure $\nu$ on the set $cX \setminus X$ with the induced topology, $\widetilde{\mu} = \mu + \nu$ is a Borel measure on $cX$ such that $\widetilde{\mu}|_X = \mu$. Since the aim of compactification is the study of a measure $\mu$ on a locally compact set $X$, we will consider only such extensions $\widetilde{\mu}$ on $cX$ that $\mu(cX \setminus X) = 0$.

\begin{lemm}\label{Smu_diff_comp}
Let $X$ be a non-compact locally compact Cantor set and $\mu \in M^0(X)$. Let $c_1X$, $c_2X$ be the compactifications of $X$ such that $c_1X \leq c_2X$. Denote by $\mu_1$ the extension of $\mu$ on $c_1 X$ and by $\mu_2$ the extension of $\mu$ on $c_2 X$.
Then $S(\mu) \subseteq S(\mu_1) \subseteq S(\mu_{2})$.
\end{lemm}

\noindent \textbf{Proof.} Since $c_1X \leq c_2X$, there exists a continuous map $f \colon c_2X \rightarrow c_1 X$ such that $f(c_2 X \setminus X) = c_1 X \setminus X$ and $fc_2(x) = c_1(x)$ for any $x \in X$. Since $f$ is continuous, it suffices to prove that $f$ preserves measure, that is $\mu_1(V) = \mu_2(f^{-1}(V))$ for any compact open $V \subset X$.
Recall that we can identify $c_i(X)$ with $X$. Hence we can consider $f$ as an identity on $X \subset c_iX$ and $f$ preserves measure.
That is, for every compact open subset $U$ of $X$ we have $\mu(U) = \mu_1(U) = \mu_2(U)$. Hence $S(\mu) \subseteq S(\mu_1)$.
Since $\mu(c_iX \setminus X) = 0$, the measure of any clopen subset of $c_iX$ is the sum of measures of compact open subsets of $X$. Hence the measures of all clopen sets are preserved. Thus, $S(\mu_1) \subseteq S(\mu_{2})$.
$\blacksquare$

\begin{remar}
We can consider the homeomorphic embedding of a set $X$ into a non-compact locally compact Cantor set $Y$ such that $\mu(Y \setminus X) = 0$. Then, by the same arguments as above, the inclusion $S(\mu|_X) \subseteq S(\mu|_Y)$ holds.
\end{remar}

\begin{thm}\label{krit_good}
Let $X$ be a non-compact locally compact Cantor set and $\mu \in M^0(X)$ be a good measure. Let $cX$ be any compactification of $X$. Then $\mu$ is good on $cX$ if and only if $S(\mu|_{cX}) \cap [0, \mu(X)) = S(\mu|_X)$.
\end{thm}

\noindent \textbf{Proof}.
First, we prove the ``if'' part. Let $V$ be a clopen set in $cX$. Consider two cases. First, let $V \cap (cX \setminus X) = \emptyset$. Then $V$ is a compact open subset of $X$. Since $\mu$ is good on $X$ and $S(\mu|_{cX}) \cap [0, \mu(X)) = S(\mu|_X)$, we see that $V$ stays good in $cX$. Now, suppose that $V \cap (cX \setminus X) \neq \emptyset$.
Then $V \cap X$ is an open set and $\mu(V) = \mu(V \cap  X) = \mu(\bigsqcup_{n=1}^{\infty} V_n)$ where each $V_n$ is a compact open set in $X$. Let $U$ be any compact open subset of $X$ with $\mu(U) < \mu(V)$. Then there exists $N \in \mathbb{N}$ such that $\mu(U) < \mu(\bigsqcup_{n = 1}^{N} V_n)$. The set $Z = \bigsqcup_{n = 1}^{N} V_n$ is a compact open subset of $X$. Since $S(\mu|_{cX}) \cap [0, \mu(X)) = S(\mu|_X)$, we have $\mu(U) \in S(\mu|_X)$. Since $\mu$ is good on $X$, there exists a compact open subset $W \subset Z$ such that $\mu(W) = \mu(U)$.

Now we prove the ``only if'' part. Assume the converse. Suppose that $\mu$ is good and the equality does not hold.
Then there exists $\gamma \in (0,\mu(X))$ such that $\gamma \in S(\mu|_{cX}) \setminus S(\mu|_X)$. Since $S(\mu|_X)$ is dense in $(0,\mu(X))$, there exists a compact open subset $U \subset X$ such that $\mu(U) > \gamma$.
Hence $\gamma \in S(\mu|_{cX}) \cap [0, \mu(U)]$ and $\gamma \not \in S(\mu|_U)$. Thus $U$ is not good and we get a contradiction. $\blacksquare$

\begin{remar}
By Proposition~\ref{basic_prop}, the set $X \setminus \mathfrak{M}_\mu$ is a non-compact locally compact Cantor set and $\overline{X \setminus \mathfrak{M}_\mu} = X$. Thus, the set $X \setminus \mathfrak{M}_\mu$ can be homeomorphically embedded into $X$ and then into some compactification $cX$. After embedding $X \setminus \mathfrak{M}_\mu$ into $X$, we add only compact open sets of infinite measure. Hence if $\mu$ was good on $X \setminus \mathfrak{M}_\mu$, it remains good on $X$ and $S(\mu|_{X \setminus \mathfrak{M}_\mu}) = S(\mu|_X)$. We can consider $X$ as a step towards compactification of $X\setminus \mathfrak{M}_\mu$ and include $\mathfrak{M}_\mu$ into $cX \setminus X$.
The measure $\mu \in M^0(X)$ is locally finite on $X \setminus \mathfrak{M}_\mu$, so we can consider only locally finite measures among infinite ones.
\end{remar}

If $\mu$ is not good on a locally compact Cantor set $X$ then clearly $\mu$ is not good on any compactification $cX$.

\begin{corol}
Let $\mu$ be a good infinite locally finite measure on a non-compact locally compact Cantor set $X$. Then $\mu$ is good on $\omega X$.
\end{corol}

\noindent \textbf{Proof.} By definition of topology on $\omega X$, the ``new'' open sets have compact complement. Since $\mu$ is locally finite on $X$, the measure of compact subsets of $X$ is finite. Hence the measure of each new clopen set is infinite. By Theorem~\ref{krit_good}, $\mu$ is good on $\omega X$.
$\blacksquare$

\begin{thm}\label{gamma}
Let $\mu$ be a good measure on a non-compact locally compact Cantor set $X$. Then for any $\gamma \in [0,\mu(X))$ there exists a compactification $cX$ such that $\gamma \in S(\mu|_{cX})$.
\end{thm}

\noindent \textbf{Proof.} The set $S(\mu|_{X})$ is dense in $[0,\mu(X))$. Hence for every $\gamma \in [0,\mu(X)$ there exist $\{\gamma_n\}_{n=1}^{\infty} \subset S(\mu|_{X})$ such that $\lim_{n\rightarrow \infty} \gamma_n = \gamma$. Since $\mu$ is good, there exist disjoint compact open subsets $\{U_n\}_{n=1}^{\infty}$ such that $\mu(U_n) = \gamma_n$. Then $U = \bigsqcup_{n=1}^{\infty} U_n$ is a non-compact locally compact Cantor set. Consider the compactification $cX = \omega U \sqcup c(X \setminus U)$, where $c(X \setminus U)$ is any compactification of $X \setminus U$.
Then $\omega U$ is a clopen set in $cX$ and $\mu(\omega U) = \gamma \in S(\mu|_{cX})$. $\blacksquare$

\medskip
From Theorems~\ref{krit_good},~\ref{gamma} the corollary follows:
\begin{corol}
For any measure $\mu$ on a non-compact locally compact Cantor space $X$ there exists a compactification $cX$ such that $\mu$ is not good on $cX$.
\end{corol}

If a measure $\mu \in M^0(X)$ is a good probability measure then, by Theorem~\ref{krit_good}, the measure $\mu$ is good on $cX$ if and only if $S(\mu|_{cX}) = S(\mu|_X) \cup \{1\}$.

\begin{prop}\label{1}
Let $X$ be a non-compact locally compact Cantor set and $\mu \in M_f(X)$. If there exists a compactification $cX$ such that $S(\mu|_{cX}) = S(\mu|_X) \cup \{1\}$ then $1 \in G(S(\mu|_X))$.
\end{prop}

\noindent \textbf{Proof}. Let $\gamma \in S(\mu|_{cX}) \cap (0,1)$. Since the complement of a clopen set is a clopen set, we have $1 - \gamma \in S(\mu|_{cX})$. Since $S(\mu|_{cX}) = S(\mu|_X) \cup \{1\}$, we have $\gamma, 1 - \gamma \in S(\mu|_{X})$. Hence $1 \in G(S(\mu|_X))$.
$\blacksquare$

\medskip
Thus, if $1 \not \in G(S(\mu|_X))$ then for any compactification $cX$ the set $S(\mu|_X)$ cannot be preserved after the extension. The examples are given in the last section.

The corollary follows from Proposition~\ref{1} and Theorem~\ref{krit_good}.
\begin{corol}
Let $\mu$ be a probability measure on a non-compact locally compact Cantor set $X$ and $1 \not \in G(S(\mu|_X))$. Then for any compactification $cX$ of $X$, $\mu$ is not good on $cX$.
\end{corol}

\begin{thm}\label{goodAlex}
Let $\mu$ be a good probability measure on a non-compact locally compact Cantor set $X$. Then $\mu$ is good on Alexandroff compactification $\omega X$ if and only if $1 \in G(S(\mu|_X))$.
\end{thm}

\noindent \textbf{Proof}. By Proposition~\ref{1} and Theorem~\ref{krit_good}, if $\mu$ is good on $\omega X$ then $1 \in G(S(\mu|_X))$.

Suppose $\mu$ is good on $X$ and $1 \in G(S(\mu|_X))$. Since $\mu$ is good, any compact open subset of $\mu$ is good, hence for every compact open $U \subset X$ we have $\mu(U) = G(S(\mu|_X)) \cap [0, \mu(U)] = S(\mu)\cap [0, \mu(U)]$. 
Every clopen subset of $\omega X$ has a compact open subset of $X$ as a complement. 
Hence for every clopen $V \subset \omega X$ we see that $\mu(V) = 1 - \mu(X \setminus V) \in G(S(\mu|_X)) \cap (0,1) = S(\mu|_X)$. So, $S(\mu|_{\omega X}) = S(\mu|_X) \cup \{1\}$. Hence $\mu$ is good on $\omega X$ by Theorem~\ref{krit_good}. $\blacksquare$

For a Cantor set $Y$ denote by $M^0(Y)$ the set of all either finite or non-defective measures on $Y$ (see~\cite{K}). Since an open dense subset of a Cantor set is a locally compact Cantor set, the corollary follows:
\begin{corol}\label{good_subs}
Let $Y$ be a (compact) Cantor set and measure $\mu \in M^0(Y)$. Let $X \subset Y$ be an open dense subset of $Y$ of full measure. If $\mu$ is good on $Y$ then $\mu$ is good on $X$.
\end{corol}

\noindent \textbf{Proof}.
The set $X$ is a locally compact Cantor set and $Y$ is a compactification of $X$. Any compact open subset $U$ of $X$ is a clopen subset of $Y$ and all clopen subsets of $U$ are compact open sets. Thus, a $\mu|_Y$-good compact open set in $X$ is, a fortiory, $\mu|_X$-good.
$\blacksquare$

\medskip
Thus, the extensions of a non-good measure are always non-good.
The corollary follows from Lemma~\ref{Smu_diff_comp}, Theorem~\ref{krit_good} and Corollary~\ref{good_subs}.

\begin{corol}\label{2compactifications}
Let $X$ be a non-compact locally compact Cantor set and $\mu \in M^0(X)$. Let $c_1X$, $c_2X$ be compactifications of $X$ such that $c_1X \geq c_2X$. Let $\mu$ be good on $c_1X$. Then $\mu$ is good on $c_2X$. Moreover, if $\mu \in M_f(X)$ then $\mu|_{c_1X}$ is homeomorphic to $\mu|_{c_2X}$.
\end{corol}

\begin{remar}
 Recall that Alexandroff compactification $\omega X$ is the smallest element in the set of all compactifications of $X$. Hence, if $\mu$ is not good on $\omega X$ then $\mu$ is not good on any compactification $cX$ of $X$.
\end{remar}

The following theorem can be proved using the results of Akin~\cite{Akin2} for measures on compact sets.
\begin{thm}
Let $X$, $Y$ be non-compact locally compact Cantor spaces, and $\mu \in M^0_f(X)$, $\nu \in M^0_f(Y)$ be good measures such that their extensions to $\omega X$, $\omega Y$ are good. Then $\mu|_X$ and $\nu|_Y$ are homeomorphic if and only if $S(\mu|_X) = S(\nu|_Y)$.
\end{thm}

\noindent \textbf{Proof.} The ``only if'' part is trivial, we prove the ``if'' part. Since $\mu|_{\omega X}$ and $\nu|_{\omega Y}$ are good by Theorem~\ref{krit_good}, we have $S(\mu|_{\omega X}) = S(\nu|_{\omega Y})$. Denote by $x_0 = \omega X \setminus X$ and $y_0 = \omega Y \setminus Y$. By Theorem 2.9~\cite{Akin2}, there exists a homeomorphism $f \colon \omega X \rightarrow \omega Y$ such that $f_*\mu = \nu$ and $f(x_0) = y_0$. Hence $f(X) = Y$ and the theorem is proved. $\blacksquare$

\medskip
In Example 1, we present a class of good measures on non-compact locally compact Cantor sets such that these measures are not good on the Alexandroff compactifications. Thus, these measures are not good on any compactification of the corresponding non-compact locally compact Cantor sets.

\section{Examples}

\begin{example}[Ergodic invariant measures on stationary Bratteli diagrams]
Let $B$ be a non-simple stationary Bratteli diagram with the matrix $A$ transpose to the incidence matrix. Let $\mu$ be an ergodic $\mathcal{R}$-invariant measure on $B$ (see~\cite{S.B.O.K.,S.B.,K}). Let $\alpha$ be the class of vertices that defines $\mu$. Then $\mu$ is good as a measure on a non-compact locally compact set $X_\alpha$. The measure $\mu$ on $X_\alpha$ can be either finite or infinite, but it is always locally finite.
The set $X_B$ is a compactification of $X_\alpha$. Since $\mu$ is ergodic, we have $\mu(X_B \setminus X_\alpha) = 0$. In~\cite{S.B.O.K.,K} the criteria of goodness for probability or non-defective measure $\mu$ on $X_B$ were proved in terms of the Perron-Frobenius eigenvalue and eigenvector of $A$ corresponding to $\mu$ (see Theorem 3.5~\cite{S.B.O.K.} for probability measures and Corollary 3.4~\cite{K} for infinite measures). It is easy to see that these criteria are particular cases of Theorem~\ref{krit_good}.

We consider now a class of stationary Bratteli diagrams and give a criterion when a measure $\mu$ from this class is good on the Alexandroff compactification $\omega X_\alpha$.
Fix an integer $N \geq 3$ and let
$$
F_N =
\begin{pmatrix}
2 & 0 & 0\\
1 & N & 1\\
1 & 1 & N \\
\end{pmatrix}
$$
be the incidence matrix of the Bratteli diagram $B_N$. For $A_N = F_N^T$ we easily find the Perron-Frobenius eigenvalue $\lambda = N+1$ and the corresponding probability eigenvector
$$
x = \left(\frac{1}{N},\ \frac{N-1}{2N},\ \frac{N-1}{2N}\right)^T.
$$
Let $\mu_N$ be the measure on $B_N$ determined by $\lambda$ and the eigenvector  $x$. The measure $\mu_N$ is good on $\omega X_\alpha$ if and only if for there exists $R \in \mathbb{N}$ such that $\frac{2(N+1)^R}{N-1}$ is an integer. This is possible if and only if $N = 2^k+1$, $k \in \mathbb{N}$. For instance, the measure $\mu_N$ is good on $\omega X_\alpha$ for $N = 3, 5$ but is not good for $N = 4$. Note that the criterion for goodness on $\omega X_\alpha$ here is the same as for goodness on $X_B$. This example is a particular case of more general result (the notation from~\cite{S.B.O.K.} is used below):
\end{example}

\begin{prop}
Let $B$ be a stationary Bratteli diagram defined by a distinguished eigenvalue $\lambda$ of the matrix $A = F^T$. Denote by $x = (x_1,...,x_n)^T$ the corresponding reduced vector. Let the vertices $2, \ldots, n$ belong to the distinguished class $\alpha$ corresponding to $\mu$. Then $\mu$ is good on $X_B$ if and only if $\mu$ is good on $\omega X_\alpha$.
\end{prop}

\noindent \textbf{Proof.}
By Theorem~\ref{Alex} and Corollary~\ref{2compactifications}, if $\mu$ is good on $X_B$ then $\mu$ is good on $\omega X_\alpha$.
We prove the converse. By Theorem 3.5 in~\cite{S.B.O.K.} and Theorem~\ref{goodAlex}, it suffices to prove that $1 \in G(S(\mu|_{X_\alpha}))$ only if there exists $R \in \mathbb{N}$ such that $\lambda^R x_1 \in H(x_2,...,x_n)$. Note that $G(S(\mu|_{X_\alpha})) = \left(\bigcup_{N=0}^\infty \frac{1}{\lambda^N} H(x_2,...,x_n) \right)$, where $H(x_2,...,x_n)$ is an additive group generated by $x_2,...,x_n$.
Suppose that $1 \in G(S(\mu|_{X_\alpha}))$. Since $\sum_{i=1}^n x_i = 1$, we see that $x_1 \in G(S(\mu|_{X_\alpha}))$, hence there exists $R \in \mathbb{N}$ such that $\lambda^R x_1 \in H(x_2,...,x_n)$. $\blacksquare$

\medskip
Return to a general case of ergodic invariant measures on stationary Bratteli diagrams.
If $\mu$ is a probability measure on $X_\alpha$ and $S(\mu|_{X_\alpha}) \cup \{1\} = S(\mu|_{X_B})$ then, by Lemma~\ref{Smu_diff_comp}, we have $S(\mu|_{\omega X_\alpha}) = S(\mu|_{X_B})$. By Theorem~\ref{krit_good}, the measure $\mu$ is good on $\omega X_\alpha$.
Hence $\mu|_{\omega X_\alpha}$ is homeomorphic to $\mu|_{X_B}$ (see~\cite{Akin2}).
If $\mu$ is infinite, then the measures $\mu|_{\omega X_\alpha}$ and $\mu|_{X_B}$ are not homeomorphic since $\mathfrak{M}_{\mu|_{\omega X_\alpha}}$ is one point and $\mathfrak{M}_{\mu|_{X_B}}$ is a Cantor set (see~\cite{K}).

\begin{example}
Let $X$ be a Cantor space and $\mu$ be a good probability measure on $X$ with $S(\mu) = \{\frac{m}{2^n} : m \in \mathbb{N} \cap [0,2^n], n \in \mathbb{N}\}$ (for example a Bernoulli measure $\beta (\frac{1}{2},\frac{1}{2})$). Clearly, $\mu_n = \frac{1}{2^n}\mu$ is a good measure for $n \in \mathbb{N}$ with $S(\mu_n) = \frac{1}{2^n} S(\mu) \subset S(\mu)$. Let $\{X_n, \mu_n\}_{n=1}^{\infty}$ be a sequence of Cantor spaces with measures $\mu_n$. Let $A = \bigsqcup_{n=1}^{\infty} X_n$ be the disjoint union of $X_n$. Denote by $\nu$ a measure on $A$ such that $\nu|_{X_n} = \mu_n$. Then $\nu$ is a good measure on a locally compact Cantor space $A$ with $S(\nu) = S(\mu) \cap [0,1)$.

Consider the one-point compactification $\omega A$ and the extension $\nu_1$ of $\nu$ to $\omega A$. We add to $S(\nu)$ the measures of sets which contain $\{\infty\}$ and have a compact open complement. Hence we add the set $\Gamma = \{1 - \gamma : \gamma \in S(\nu)\}$. Since $\Gamma \subset S(\nu) \cup \{1\}$, we have $S(\nu_1) = S(\nu) \cup \{1\}$. By Theorem~\ref{krit_good}, the measure $\nu_1$ is good on $\omega A$.

Consider the two-point compactification of $A$. Let $A = A_1 \sqcup A_2$ where $A_1 = \bigsqcup_{k = 1}^{\infty} X_{2k - 1}$ and $A_2 = \bigsqcup_{j = 1}^{\infty} X_{2j}$. Then $cA = \omega A_1 \sqcup \omega A_2$ is a two-point compactification of $A$. Let $\nu_2$ be the extension of $\nu$ to $cA$. Then $\nu_2(A_1) = \frac{2}{3} \not \in S(\nu)$. Hence, by Theorem~\ref{krit_good}, the measure $\nu_2$ is not good on $cA$.

In the same example, we can make a two-point compactification which preserves $S(\nu|_A)$. Since $\mu_n$ is good for $n \in \mathbb{N}$, there is a compact open partition $X_n^{(1)} \sqcup X_n^{(2)} = X_n$ such that $\mu_n (X_n^{(i)}) = \frac{1}{2^{n+1}}$ for $i = 1,2$. Let $B_i = \bigsqcup_{n=1}^{\infty} X_n^{(i)}$ for $i = 1,2$. Consider $\widetilde{c}A = \omega B_1 \sqcup \omega B_2$. Then it can be proved the same way as above that $S(\nu|_{\widetilde{c}A}) = S(\nu|_A) \cup \{1\}$.
\end{example}

\begin{example}
Let $\mu = \beta(\frac{1}{3}, \frac{2}{3})$ be a Bernoulli (product) measure on Cantor space $Y = \{0,1\}^\mathbb{N}$ generated by the initial distribution $p(0) = \frac{1}{3}$, $p(1) = \frac{2}{3}$. Then $\mu$ is not good but $S(\mu) = \{\frac{a}{3^n} : a \in \mathbb{N} \cap [0, 3^n], n \in \mathbb{N}\}$ is group-like (see~\cite{Akin1}).
Let $X$ be any open dense subset of $Y$ such that $\mu(Y \setminus X) = 0$. Thus, $Y$ is a compactification of a non-compact locally compact Cantor space $X$ and $\mu$ extends from $X$ to $Y$. Then $\mu$ is not good on $X$.

The compact open subsets of $X$ are exactly the clopen subsets of $Y$ that lie in $X$. The compact open subset of $X$ is a union of the finite number of compact open cylinders. Consider any compact open cylinder $U = \{a_0,...,a_n,*\}$ which consists of all points in $z \in Y$ such that $z_i = a_i$ for $0 \leq i \leq n$. Then $U$ is a disjoint union of two subcylinders $V_1 = \{a_0,...,a_n,0,*\}$ and $V_2 = \{a_0,...,a_n,1,*\}$ with $\mu(V_2) = 2 \mu(V_1)$.
Let the numerator of the fraction $\mu(V_1)$ be $2^k$. Then the numerator of the fraction $\mu(V_2)$ is $2^{k+1}$. Moreover, for any compact open $W \subset V_2$ the numerator of the fraction $\mu(W)$ will be divisible by $2^{k+1}$. Since $S(\mu)$ contains only finite sums of measures of cylinder compact open sets and the denominators of elements of $S(\mu)$ are the powers of $3$, there is no compact open subset $W \subset V_2$ such that $\mu(W) = \mu(V_1)$. Hence $\mu$ is not good on $X$.

Moreover, let $x = \{00...\}$ be a point in $Y$ which consists only of zeroes.
Then $S(\mu |_{Y \setminus\{x\}})\varsubsetneq S(\mu |_Y)$ while $S(\mu |_{Y \setminus\{y\}}) = S(\mu |_Y)$ for any $y \neq x$.

Consider the case $y \neq x$.
Let, for instance, $y = \{111....\}$, all other cases are proved in the same way.
Let $U_n = \{z \in Y : z_0 = ... = z_{n-1} = 1, z_{n} = 0\}$. Then $Y \setminus\{y\} = \bigsqcup_{n=1}^\infty U_n \sqcup \{0*\}$. Denote by $S_N = \mu(\bigsqcup_{n=1}^N U_n)$ and $S_0 = 0$. Then $\lim_{N \rightarrow \infty} S_N = \frac{2}{3}$. Let $G = \{\frac{a}{3^n} : a \in \mathbb{Z}, n \in \mathbb{N}\}$. Then $G$ is an additive subgroup of reals and $S(\mu|_Y) = G \cap [0,1]$. We prove that $S(\mu |_{Y \setminus\{y\}}) = G \cap [0,1)$, i.e. for every $n \in \mathbb{N}$ and $a = 0,..., 3^n - 1$ there exists a compact open set $W$ in $Y \setminus \{y\}$ such that $\mu(W) = \frac{a}{3^n}$. Indeed, we have $S(\mu|_{\{0*\}}) = G \cap [0, \frac{1}{3}]$ and $[0,1) = \cup_{n=0}^\infty [S_N, S_N + \frac{1}{3}]$. Hence $G \cap [0,1) = \cup_{n=0}^\infty (G \cap [S_N, S_N + \frac{1}{3}])$. For every $\gamma \in G$ there exists $N\in \mathbb{N}$ such that $\gamma \in [S_N, S_N + \frac{1}{3}]$. There exists a compact open subset $W_0$ of $\{0*\}$ such that $\mu(W_0) = \gamma - S_N$. Set $W = U_N \sqcup W_0$. Then $W$ is a compact open subset of $Y \setminus \{y\}$ and $\mu(W) = \gamma$.

Now consider the set $Y \setminus\{x\}$. Every cylinder that lies in $Y \setminus\{x\}$ has even numerator, hence $S(\mu |_{Y \setminus\{x\}})\varsubsetneq S(\mu |_Y)$. It can be proved in the same way as above that $S(\mu |_{Y \setminus\{x\}}) = \{\frac{2k}{3^n} : k \in \mathbb{N}\} \cap [0,1)$.

\end{example}

\begin{example}[$(C,F)$-construction]
Denote by $|A|$ the cardinality of a set $A$.
Given two subsets $E, F \subset \mathbb{Z}$, by $E+F$ we mean $\{e+f | e \in E, f \in F\}$ (for more details see~\cite{D1,D2}).
Let $\{F_n\}_{n=1}^\infty, \{C_n\}_{n=1}^\infty \subset \mathbb{Z}$ such that for each $n$

(1) $|F_n|<\infty$, $|C_n|<\infty$,

(2) $|C_n| > 1$,

(3) $F_n + C_n + \{-1,0,1\} \subset F_{n+1}$,

(4) $(F_n + c) \cap (F_n + c') = \emptyset$ for all $c \neq c' \in C_{n+1}$.

Set $X_n = F_n \times \prod_{k > n} C_k$ and endow $X_n$ with a product topology. By (1),(2), each $X_n$ is a Cantor space.

For each $n$, define a map $i_{n,n+1} \colon X_{n} \rightarrow X_{n + 1}$ such that
$$
i_{n,n+1}(f_n, c_{n+1}, c_{n+2},...) = (f_n + c_{n+1}, c_{n+2},...).
$$
By (1), (2) each $i_{n,n+1}$ is a well defined injective continuous map.
Since $X_n$ is compact, we see that $i_{n,n+1}$ is a homeomorphism between $X_n$ and $i_{n,n+1}(X_n)$. So the embedding $i_{n,n+1}$ preserves topology.
The set $i_{n, n+1}(X_n)$ is a clopen subset of $X_{n+1}$. Let $i_{m,n} \colon X_m \rightarrow X_n$ such that  $i_{m,n} = i_{n, n-1} \circ i_{n-1, n-2} \circ ... \circ i_{m+1,m}$ for $m < n$ and $i_{n,n} = id$.
Denote by $X$ the topological inductive limit of the sequence $(X_n, i_{n,n+1})$. Then $X = \bigcup_{n=1}^\infty X_n$.
Since $i_{m,n} = i_{n, n-1} \circ i_{n-1, n-2} \circ ... \circ i_{m+1,m}$ for $m < n$, we can write $X_1 \subset X_2 \subset ...$
The set $X$ is a non-compact locally compact Cantor set.
The Borel $\sigma$-algebra on $X$ is generated by cylinder sets $[A]_n = \{x \in X : x = (f_n, c_{n+1},c_{n+2},...) \in X_n \mbox{ and } f_n \in A\}$. There exists a canonical measure on $X$. Let $\kappa_n$ stand for the equidistribution on $C_n$ and let $\nu_n = \frac{|F_n|}{|C_1|...|C_n|}$ on $F_n$. The product measure on $X_n$ is defined as $\mu_n = \nu_n \times \kappa_{n+1} \times \kappa_{n+2}\times ...$ and a $\sigma$-finite measure $\mu$ on $X$ is defined by restrictions $\mu|_{X_n} = \mu_n$. The measure $\mu$ is a unique up to scaling ergodic locally finite invariant measure for a minimal self-homeomorphism of $X$ (for more details see~\cite{D1, D2}). For every two compact open subsets $U,V \subset X$ there exists $n \in \mathbb{N}$ such that $U, V \subset X_n$.
The measure $\mu$ is obviously good, since the restriction of $\mu$ onto $X_n$ is just infinite product of equidistributed measures on $F_n$ and $C_m$, $m > n$. We have $S(\mu) = \{\frac{a}{|C_1|...|C_n|} : a, n \in \mathbb{N}\} \cap [0, \mu(X))$.

\end{example}

\begin{example}
Let $p$ be a prime number and $\mathbb{Q}_p$ be the set of $p$-adic numbers. Endowed with the $p$-adic norm, the set $\mathbb{Q}_p$ is a non-compact locally compact Cantor space. Then the Haar measure $\mu$ on $\mathbb{Q}_p$ is good and $S(\mu) = \{n p^{\gamma} | n \in \mathbb{N}, \gamma \in \mathbb{Z}\}$.
\end{example}

\Large \noindent \textbf{Acknowledgement}

\medskip
\normalsize
I am grateful to my advisor Sergey Bezuglyi for giving me the idea of this work, for many helpful discussions and for reading the preliminary versions of this paper.

\end{document}